\documentclass{amsart}
\theoremstyle{plain}
\newtheorem{thm}{Theorem}

\errorcontextlines=0 \numberwithin{equation}{section}

\begin{document}
\large
\title[Mass and Flows]
{When is the Hawking mass monotone under Geometric Flows}
\author{John Bland and Li Ma }

\address{JB: Department of Mathematics\\
Toronto University\\
Toronto, Ontario, Canada M5S 2E4}

\address{LM: Department of Mathematical Sciences \\
Tsinghua University \\
Beijing 100084 \\
China}

\email{lma@math.tsinghua.edu.cn}

\dedicatory{}
\date{}
\thanks{The research is partially supported by the National Natural Science
Foundation of China 10631020 and SRFDP 20060003002}


\begin{abstract}
In this paper, we study the relation of the monotonicity of Hawking
Mass and geometric flow problems. We show that along the
Hamilton-DeTurck flow with bounded curvature coupled with the
modified mean curvature flow, the Hawking mass of the hypersphere
with a sufficiently large radius in Schwarzschild spaces is monotone
non-decreasing.

\end{abstract}
\maketitle

\section{Introduction}\label{ONE}

In this work we plan to study the monotonicity of Hawking mass under
the Ricci flow coupled with the geometric flow for hypersurfaces in
a slice of the Ricci flow. Actually, H.Bray \cite{Bray} asked if
there exists a flow which makes the quasi-local mass functional
nondecreasing.

Given a time interval $I=(a,b)$. The Ricci flow on a manifold
$M^{n+1}$ is defined by a one parameter family of metrics $\{g(t);
t\in I\}$ satisfying
\begin{equation}\label{rf}
\partial_tg(t)=-2Rc(g(t)),\; \; t\in I
\end{equation}
where $Rc(g(t))$ is the Ricci tensor of the metric $g(t)$. This flow
was introduced by R.Hamilton in 1982 (\cite{H95}, \cite{M}) and was
used by G.Perelman (\cite{P}) in 2002-3 to give an outline of the
proof of Poincare conjecture in dimension three. Ricci flow is a
fully nonlinear equation with beautiful geometrical applications
(see \cite{H95}, \cite{P}, also \cite{CM05}). The local existence
result on a compact manifold was obtained by R.Hamilton in the
famous paper in 1982. In the paper \cite{DM1}, we have proved that
the asymptotically locally Euclidean (ALE) property is preserved by
the Ricci flow (and the related Hamilton-Ricci flow, which is
equivalent to the Ricci flow flow after a change of gauge). An
interesting convergence result is obtained in \cite{OW}. The Ricci
flow was also been discussed in Physical literature, such as in
quantum filed theory, where it was considered as an approximation to
the renormalization group (RG) flow for the 2-dimensional nonlinear
sigma model.

When interpreted as the statistical physics analogy to Einstein
filed equation, Ricci flow can be considered as the gradient flow
of the W-functional, which is a non-decreasing entropy functional
along the Ricci flow coupling a heat equation. For a precise
definitions of W-functional, one may look at G.Perelman's paper
\cite{P}.

We want to find such a flow (which will be the modified mean
curvature flow, see \ref{NMCF} in section \ref{THREE} for the
definition and see\cite{HI01} for related flow) so that the Hawking
mass (\cite{HY86}, \cite{QT07}, \cite{WY}, and\cite{W02}) is
monotone along it. Given a Riemannian manifold $(M^{n+1},g)$. Assume
that $n=2$. Let's recall the definition of Hawking mass for a
hypersurface $F: \Sigma \to M$. Let $\nu$ be the outward unit normal
vector to $\Sigma$. Let $H=H(\Sigma)$ and $V=V(\Sigma)$ be the mean
curvature of $\Sigma $and the area of $\Sigma$ in $M$ respectively.
Then the Hawking mass of the hypersurface $\Sigma$ is defined by
$$
m_{H}=m(\Sigma)=\sqrt{\frac{{V}}{16\pi}}(1-\frac{1}{16\pi}\int_{\Sigma}|H|^2).
$$

Given a family of hypersurfaces $F(\cdot, t)$. We want to find a
real valued smooth function $p(F)$ such that along the geometric
flow of the form
$$
\partial_t F=p(F)\nu(F),
$$
the Hawking mass is monotone. This is a natural question. Usually,
one should have a lot such kind of flows. We just pick up a natural
one. Our main result is Theorem \ref{Hawkingone} in section
\ref{THREE}. Roughly speaking, the result says that along the
Hamilton-DeTurck flow with bounded curvature coupled with the
modified mean curvature flow, the Hawking mass of the hypersphere
with a sufficiently large radius in Schwarzschild space is monotone
non-decreasing. We should say that our result may be extended to
non-symmetric cases, however, our aim here is to introduce the flow.
This result is interesting in general relativity.

The plan of this paper is the following. In the next section, we
give useful geometric formulae  and properties of the Ricci flow for
warped product spaces. In section \ref{THREE}, we prove the main
result stated in Theorem \ref{Hawkingone}.

\section{The flow on warped product spaces and related}\label{TWO}

  We denote
by $R(t)$ the scalar curvature of the metric $g(t)$. Then the
function $R(t)$ satisfies the following evolution equation
\begin{equation}\label{Scflow}
\partial_t R=\Delta R+2|Rc|^2\geq \Delta R+\frac{2}{n+1}R^2.
\end{equation}
This is an important differential inequality which give us the blow
up in finite time in compact manifold with initial positive scalar
curvature. In particular, we have that the positivity of scalar
curvature is preserved along the Ricci flow.
However, we assume that
the Ricci flow has uniformly bounded curvature (see \cite{DM1} and
\cite{Shi89}).

 Let
$M=\mathbb{R}\times S^n$ with the product metric
\begin{equation}\label{model1}
g=\phi(x)^2dx^2+\psi(x)^2\hat{g}.
\end{equation}
Here $S^n$ be the standard n-sphere so that $\hat{g}=d\sigma^2$.
We shall write by $x=r$ as the standard radial variable in the
Schwarzschild coordinates.

Introduce the arc-length parameter $s$ and sectional curvatures in
the following way. Let
\begin{equation}\label{length}
s=\int_0^x\phi(\tau)d\tau.
\end{equation}
Let $$ K_0=-\frac{\psi_{ss}}{\psi}.
$$ be the sectional curvature of the 2-planes perpendicular to the
spheres $\{x\}\times S^n$
 and
$$
K_1=\frac{1-\psi_s^2}{\psi^2}
$$
the sectional curvature of the 2-planes tangential to the spheres.
Note that
$$
\frac{\partial}{\partial
s}=\frac{1}{\phi(x)}\frac{\partial}{\partial x}.
$$

In terms of the arc-length parameter $s$, the metric can be read
as
$$ g=ds^2+\psi^2\hat{g}.
$$

 It is known
that the Ricci curvature tensor of the product metric $g$ is given
by
$$
Rc=n \{-\frac{\psi_{xx}\psi+\psi_x\phi_x}{\psi\phi}
\}dx^2+\{-\frac{\phi\psi_{xx}-(n-1)\phi\psi_x^2+\psi\phi_x\psi_x}{\phi^3}+n-1\}\hat{g}.
$$
See \cite{AK} and \cite{SM}.

 Let
$$ g_s=\psi^2\hat{g}.
$$
Then we have the hessian formula for $s$ (see \cite{PP}):
$$
Hess s=\psi \psi_s \hat{g}=\frac{\psi_s}{\psi}g_s,
$$
and the Ricci tensor can be written as
$$
Rc=nK_0ds^2+[K_0+(n-1)K_1]\psi^2\hat{g}.
$$
See \cite{SM}.

Note that the scalar curvature is given by
$$
R=nK_0+n[K_0+(n-1)K_1].
$$

Then for the warped product metrics, the Ricci flow (\ref{rf})
becomes the system
\begin{equation}\label{Rcst}
\left\{\begin{array}{ll}
\psi_t=-[K_0+(n-1)K_1]\psi=\psi_{ss}-(n-1)\frac{1-\psi_s^2}{\psi},
 \\
\phi_t =-nK_0\phi=n\frac{\psi_{ss}}{\psi}\phi.
\end{array}
\right.
\end{equation}
In the xt coordinates, we have
\begin{equation}\label{Rcxt}
\left\{\begin{array}{ll} \psi_t=
\frac{1}{\phi(x)}\frac{\partial}{\partial
x}\frac{1}{\phi(x)}\frac{\partial}{\partial x}\psi
-(n-1)\frac{1-(\frac{1}{\phi(x)}\frac{\partial}{\partial
x}\psi)^2}{\psi},
 \\
\phi_t = n\frac{\frac{\partial}{\partial
x}\frac{1}{\phi(x)}\frac{\partial}{\partial x}\psi}{\psi}.
\end{array}
\right.
\end{equation}

We now discuss some related geometric quantities. For a fixed
$r>0$, let
$$
\Sigma:=\partial B(r)=\{p\in M; s(p)=r\}
$$
be the sphere of radius $r$ in the warped Riemannian manifold
$(M,g)$. Note that the volume of $\partial B(r)$ is
$$
V(r)=|\Sigma|=\int_{\Sigma}d\sigma=\int_N\psi(s)^{n}dv_{\hat{g}}=\psi(r)^nV(N).
$$
Here $V(N)$ is the volume of $N=S^n$ in the metric $\hat{g}$. In
particular, for $n=2$, we have $V(N)=4\pi$.

 Let
$$ g_s=\psi^2\hat{g}.
$$
Then we have the hessian formula for $s$ (see \cite{PP}):
$$
Hess (s)=\psi \psi_s \hat{g}=\frac{\psi_s}{\psi}g_s.
$$
 By the hessian formula for the function $s$, we have the mean
curvature of the sphere $\Sigma$ given by
\begin{equation}\label{hessian}
H=\Delta s=n\frac{\psi_s}{\psi}=\frac{V'(s)}{V(s)}.
\end{equation}
Note that the Hawking mass of $\Sigma$ ($n=2$) is
$$
m=m(\Sigma)=\sqrt{\frac{V(s)}{16\pi}}(1-\frac{1}{16\pi}\frac{V'(s)^2}{V(s)}),
$$
which can also be written as
\begin{equation}\label{hawking}
m=\frac{\psi(s)}{2}(1-\psi^{'2}).
\end{equation}

In the following, we always assume that $n=2$.

\section{Monotonicity in Schwarzschild spaces}\label{THREE}

Recall that a Schwarzschild space has the metric of the form
\begin{equation}\label{model2}
g=exp(2\lambda(r))dr^2+r^2(d\theta^2+\sin^2\theta d\phi^2),
\end{equation}
which is assumed to be asymptotically locally Euclidean (in short,
ALE).

To keep the specific form of the metrics, we need to use the
so-called Hamilton-DeTurck flow. The Hamilton-DeTurck flow is the
following
$$
\partial_tg=-2Rc(g)+L_Xg
$$
where $X$ is the vector-field on the manifold $M$ and $L_Xg$ is
the Lie-derivative of the metric $g(t)$. This flow is important in
the sense that we can keep the form of the metrics (\ref{model2})
along it by suitable vector field $X$. In fact, we choose
$$X=[\frac{n-1}{r}(exp(2\lambda(r,t))-1)+\partial_r\lambda]\partial_r.$$

Given a manifold $M\subset R^n$ with the standard Euclidean polar
coordinates $(r,\theta)$. For the Hamilton-DeTurck Ricci flow of
the form
\begin{equation}\label{schwarz}
g(t)=exp(2\lambda(r,t))dr^2+r^2(d\theta^2+\sin^2\theta d\phi^2),
\end{equation}
we have
\begin{equation}\label{FE}
\partial_tf =\frac{1}{f^2}\frac{\partial^2 f}{\partial
r^2}-\frac{2}{f^3}(f_r)^2+(\frac{n-1}{r}
-\frac{1}{rf^2})f_r-\frac{(n-1)}{r^2f}(f^2-1),
\end{equation}
for $f=exp(\lambda)$ with $f(0)=1$ and $f(\infty)=1$. Using the
Maximum principle trick, we have that $f$ is uniformly bound both in
time variable and space variables, which implies the global
existence of the flow (see also \cite{OW}). Invoking the ALE
properties of the Ricci flow \cite{DM1} and the gradient estimate of
Shi \cite{Shi89}, we know that
$$ f\to 1, \; \; r\partial_rf\to 0, \; \; r^2\partial_{r}^2f\to 0
$$
uniformly for $r\to\infty$; in particular,the derivative $f_t$ is
also uniformly bounded.

 Note that
the unit normal $\nu$ to $\Sigma$ is
$$
\nu=\frac{1}{\phi}\frac{\partial}{\partial x}.
$$
Hence the geometric flow for hypersurfaces can be written in the
following form:
\begin{equation}\label{Malige2}
\partial_ts=p(H(s)).
\end{equation}
Here we have
$$
H(s)=\frac{2}{rf}.
$$

In the state case when the metric is of the form (\ref{schwarz}),
we have
$$
\psi(s)=r(s), \; \; \; \frac{\partial r}{\partial
s}=e^{-\lambda}=f^{-1}.
$$
Then we have
$$
K_0=-\frac{f_r}{rf^3}, \; \; K_1=\frac{1-f^{-2}}{r^2},
$$
and
$$
R=4K_0+2K_1=-\frac{4f_r}{rf^3}+2\frac{1-f^{-2}}{r^2}.
$$

We also have $$ m(r)=\frac{r}{2}(1-f^{-2}),
$$ and
$$
\frac{d}{dr}m(r) =\frac{1-f^{-2}}{2}-rf^{-3}f_r.
$$
We remark that
$$
m(r)=\frac{r^3K_1}{2}
$$
for this metric.

In the evolving case, in some sense, we need to understand the
evolution of the curvature $K_1$. Here, we want to consider in
another way and we shall determine a geometric nature choice of
$p(r)$ to make the Hawking mass monotone.

Note that along the Hamilton-DeTurck flow coupled with the
geometric flow (\ref{Malige2})
$$
\frac{d m(r,t)}{dt}=\frac{\partial m}{\partial
r}\frac{dr}{ds}\frac{ds}{dt}+\partial_tm
=[\frac{1-f^{-2}}{2}-rf^{-3}f_r]f^{-1}p(H)+\partial_tm.
$$
Note that
$$
\partial_tm
=rf^{-3}f_t.
$$

So, putting all these together, we have
$$
\frac{d
m(r,t)}{dt}=[\frac{1-f^{-2}}{2}-rf^{-3}f_r]f^{-1}p(H)+rf^{-3}f_t.
$$
Then we have
$$
\frac{dm(r,t)}{dt}=\frac{1}{4}r^2Rf^{-1}p(H)+rf^{-3}f_t.
$$

Assume that
\begin{equation}\label{NMCF}
p(H)=R^{-1}H,
\end{equation}
and we shall call such a geometric flow as \emph{the negative mean
curvature flow}. Using $fH=\frac{2}{r}$, we have
$$
\frac{d m(r,t)}{dt} =\frac{r}{2}+rf^{-3}f_t+0(1).
$$
Note that $\frac{1}{H}=\frac{r}{2}+0(1)$ is the leading term in the
expression above. Then
$$
\frac{d m(r,t)}{dt}>0
$$
for large radius $r>0$. Therefore, we get the following result.

\begin{thm}\label{Hawkingone} Assume that the scalar curvature
$R(g_0)$ is positive for the initial Schwarzschild metric
$$
g_0=exp(2\lambda(r))dr^2+r^2(d\theta^2+\sin^2\theta d\phi^2),
$$
which is assumed to be asymptotically locally Euclidean.
 Along the Hamilton-DeTurck flow with bounded curvature coupled with the
modified Mean curvature flow, there exists a large radius $r_0>>1$
depending on the bound of curvature of Hamilton-DeTurck flow such
that for all $r>r_0$, the Hawking mass is monotone non-decreasing.
\end{thm}

It is quite clear that we can make the Hawking mass monotone
non-increasing if we choose $$ p(H)=-R^{-1}H.
$$
However, in this case, the geometric flow is a back-ward flow.
 It is an interesting question to study the same problem along the
Cross-curvature flow introduced by Chow-Hamilton \cite{MC}. We
believe the similar result is also true for general asymptotically
flat warped product spaces. The progress is still in the future.

\bigskip

 {\bf Acknowledgements}.  This work was done while the second named
author was visiting the Mathematical Department of Toronto
University, Canada, for which he would like to thank the
hospitality.

\end{document}